\documentclass[10pt,reqno]{amsart}

\usepackage[latin1]{inputenc}
\usepackage{amssymb,amsmath,amscd,amsfonts,color, amsthm}
\title[Limiting distribution of stationary Gram matrices]{The Empirical Eigenvalue Distribution of a 
Gram Matrix: From Independence to Stationarity}

\author[Hachem et al.]{W. Hachem, P. Loubaton and J. Najim}
\date{\today}

\newtheorem{theo}{Theorem}[section]

\newtheorem{lemma}[theo]{Lemma}

\newtheorem{prop}[theo]{Proposition}
\newtheorem{assump}{Assumption A-\hspace{-0.15cm}}
\newcommand{\ii}{\mathbf{i}}
\newcommand{\indicej}{{\tiny \begin{array}{l}
j_1=0:N-1\\
j_2=0:n-1
\end{array}}}

\newcommand{\EE}{\mathbb{E}}
\newcommand{\ZZ}{\mathbb{Z}}
\newcommand{\bdm}{\begin{displaymath}}
\newcommand{\edm}{\end{displaymath}}
\newcommand{\bea}{\begin{eqnarray*}}
\newcommand{\eea}{\end{eqnarray*}}

\newcommand{\Cplus}{\mathbb{C}^+}
\newcommand{\C}{\mathbb{C}}
\newcommand{\pitilde}{\tilde{\pi}}

\newcommand{\im}{\mathrm{Im}}

\newcommand{\T}{\mathrm{T}}

\numberwithin{equation}{section}
\theoremstyle{remark}
\newtheorem{rem}{Remark}[section]

\begin{document}
\bibliographystyle{plain}
\begin{abstract}
Consider a $N\times n$ matrix $Z_n=(Z^n_{j_1 j_2})$ where the individual entries are a realization of a 
properly rescaled stationary gaussian random field:
$$
Z^n_{j_1 j_2}=\frac{1}{\sqrt{n}}\sum_{(k_1,k_2)\in \mathbb{Z}^2} h(k_1,k_2)U(j_1-k_1,j_2-k_2),
$$
where $h\in\ell^1(\mathbb{Z}^2)$ is a deterministic complex summable sequence and 
$(U(j_1,j_2); (j_1,j_2)\in~\mathbb{Z}^2)$ is a sequence of independent 
complex gaussian random variables with mean zero and unit variance. 

The purpose of this article is to study the limiting empirical distribution of the eigenvalues of Gram random matrices such as 
$Z_n Z_n ^*$ and $(Z_n +A_n)(Z_n +A_n )^*$ where $A_n$ is a deterministic matrix with appropriate assumptions in the case 
where $n\rightarrow \infty$ and $\frac Nn \rightarrow c \in (0,\infty)$. 

The proof relies on related results for matrices with independent but not identically distributed entries and substantially
differs from related works in the literature (Boutet de Monvel et al. \cite{BKV96}, Girko \cite{Gir01a}, etc.).
\end{abstract}

\maketitle

\section{Introduction} 
\subsubsection*{The model} 
Let $Z_n=(Z^n_{j_1 j_2}, 0\le j_1 < N,\ 0 \le j_2 < n)$ be a 
$N \times n$ random matrix with entries 
$$
Z^n_{j_1 j_2}=\frac{1}{\sqrt{n}}\sum_{(k_1,k_2)\in \mathbb{Z}^2} h(k_1,k_2)U(j_1-k_1,j_2-k_2),
$$ 
where $(U(j_1,j_2),\ (j_1,j_2)\in \mathbb{Z}^2)$ is a sequence 
of independent complex Gaussian random variables (r.v.) such that
$\EE U(j_1,j_2) = 0$, $\EE U(j_1,j_2)^2 = 0$ and 
$\EE \left| U(j_1,j_2) \right|^2 = 1$, and  
$(h(k_1,k_2),\ (k_1,k_2)\in \ZZ^2)$ is a deterministic complex 
sequence satisfying 
$$
\sum_{(k_1,k_2)\in \mathbb{Z}^2} |h(k_1,k_2)| <\infty \ . 
$$
The bidimensional process $Z^n_{j_1 j_2}$ is a stationary gaussian field. Indeed, 
$\mathrm{cov}(Z^n_{j_1 j_2},Z^n_{j'_1 j'_2})=n^{-1}C(j_1-j'_1,j_2-j'_2)$ where
\begin{equation}\label{variance}
C(j_1,j_2)=\sum_{(k_1,k_2)\in \mathbb{Z}^2} h(k_1,k_2)h^*(k_1-j_1,k_2-j_2) 
\end{equation}
(we denote by $a^*$ the complex conjugate of $a \in \mathbb{C}$ - we also 
denote by $A^*$ the hermitian adjoint of matrix $A$).

\subsubsection*{The main results} The purpose of this article is to establish the convergence of the 
empirical distribution of the eigenvalues of various Gram matrices
based on $Z_n$. More precisely, we shall study the convergence of the
spectral distribution of $Z_n Z_n^*$ and $(Z_n+A_n)(Z_n +A_n)^*$
where $A_n$ is a deterministic matrix with a given structure. In particular, if $Z_n$ is square, we 
take $A_n$ to be Toeplitz. The contribution of this article is to provide 
a new method to study Gram matrices based on Gaussian fields. The main idea is to approximate 
the matrix $Z_n$ by a matrix $\tilde{Z}_n$ unitarily congruent to a matrix with independent but not identically distributed
entries. This method will allow us to revisit the centered case $Z_n Z_n^*$, already studied by Boutet de Monvel et al. 
in \cite{BKV96} and to establish the limiting spectral distribution of the non-centered case $(Z_n +A_n)(Z_n+A_n)^*$ 
for some deterministic matrix $A_n$.

\subsubsection*{Motivations} The motivations for such a work are twofold. First of all, we believe that this line of proof
is new. 
Let us briefly describe the three main elements of it.\\
\indent The first one is a periodization scheme popular in signal
processing and described as follows:
$$
\tilde{Z}_n =(\tilde{Z}^n_{j_1 j_2})\quad \textrm{where} \quad \tilde{Z}^n_{j_1 j_2}=\frac {1}{\sqrt{n}}\sum_{(k_1,k_2)\in 
\mathbb{Z}^2} h(k_1,k_2)
U\left( (j_1-k_1) \  \mathrm{mod}\ N, (j_2-k_2)\ \mathrm{mod}\ n\right), 
$$ 
where $\mathrm{mod}$ denotes modulo. \\
\indent The second element is an inequality due to Bai
\cite{Bai93b} involving the L\'evy distance ${\mathcal L}$
between distribution functions:
$$
{\mathcal L}^4(F^{A A^*}, F^{B B^*}) \le \frac{2}{N^2} \mathrm{Tr}(A-B)(A-B)^* \mathrm{Tr}(A A^* + B B^*),
$$
where $F^{A A^*}$ denotes the empirical distribution function of
the eigenvalues of the matrix $A A^*$ and $\mathrm{Tr}(X)$ denotes the
trace of matrix $X$. With the help of this inequality, we shall prove
that $Z_n$ and $\tilde{Z}_n$ have the same
limiting spectral distribution.\\
\indent The third element comes from the advantage of considering
$\tilde{Z}_n$. In fact, $\tilde{Z}_n$ is congruent (via Fourier
unitary transforms) to a random matrix with independent but not
identically distributed entries.  Therefore, we can (and will) rely on
results established in \cite{HLN04pre} for Gram matrices with independent but not
identically distributed
entries.\\
\indent The second motivation comes from the field of wireless
communications.  In a communication system employing antenna arrays at
the transmitter and at the receiver sides, random matrices extracted
from Gaussian fields are often good models for representing the radio
communication channel.  In this course, the stationary model as
considered above is often a realistic channel model.  The computations
of popular receiver performance indexes such as Signal to Interference
plus Noise Ratio or Shannon channel capacity heavily rely on the
knowledge of the limiting spectral distribution of matrices of the
type $Z_n Z_n^*$ (see \cite{Chu02},\cite{Liu03} and also the tutorial
\cite{TulVer04} for further references).

\subsubsection*{About the literature} Various Gram matrices based on Gaussian fields have already been studied in the literature.
The study of the general case $(Z_n+A_n)(Z_n+A_n)^*$ has been
undertaken by Girko in \cite{Gir01a}. Since no assumptions are done on
the structure of $A_n$, there might not be any limiting spectral
distribution. Girko finds asymptotic approximations of the Stieltjes
transform of $(Z_n+A_n)(Z_n+A_n)^*$. The method developed in
\cite{Gir01a} is based on an exhaustive study of each entry of the resolvent
$((Z_n+A_n)(Z_n+A_n)^* -zI)^{-1}$ added to the property that 
sufficiently remote entries are asymptotically independent.\\
\indent Boutet de Monvel et al. \cite{BKV96} have also studied Gram matrices based on stationary Gaussian fields in the case where 
the matrix has the form $V_n +Z_n Z_n^*$, $V_n$ being a 
deterministic Toeplitz matrix. Their line of proof is based on a direct study
of the resolvent, taking advantage of the gaussianity of the entries.  

\subsubsection*{Disclaimer} In this paper, we study in detail the case 
where the entries of matrix $Z_n$ are complex.
In the real case, the general framework of the proof works as well 
if one considers the real counterpart of the Fourier unitary transforms, however the computations are more involved.
We provide some details in Section \ref{section:real}.\\

\section{Assumptions and useful results}\label{resultats}

\subsection{Notations, Assumptions, Stieltjes transforms and Stieltjes kernels}
\label{sub-notations}
Let $N=N(n)$ be a sequence of integers such that 
$$
\lim_{n\rightarrow \infty} \frac{N(n)}{n} =c.
$$
We denote by $\ii$ the complex number $\sqrt{-1}$, 
by $\mathbf{1}_A(x)$ the indicator function 
over set $A$ and by $\delta_{x_0}(x)$ the Dirac measure at point $x_0$. 
A sum will be equivalently written as $\sum_{k=1}^n$ or $\sum_{k=1:n}$. 
We denote by $\mathcal{CN}(0,1)$ 
the distribution of the Gaussian complex random variable $U$ satisfying
$\EE U = 0$, $\EE U^2 = 0$, and $\EE \left| U \right|^2 = 1$ (equivalently,  
$U=A+\ii B$ where $A$ and $B$ are real independent Gaussian r.v.'s 
 with mean $0$ and standard deviation $\frac{1}{\sqrt{2}}$ each).
\begin{assump}\label{hypo1}
The entries $(Z_{j_1 j_2}^n\ ,\ 0\le j_1 < N,\,0\le j_2 < n\,,\, n\ge1)$ 
of the $N\times n$ matrix $Z_n$ are random variables defined as:
$$
Z^n_{j_1 j_2}=\frac{1}{\sqrt{n}}\sum_{(k_1,k_2)\in \mathbb{Z}^2} h(k_1,k_2)U(j_1-k_1,j_2-k_2),
$$
where $(h(k_1,k_2),\ (k_1,k_2)\in \mathbb{Z}^2)$ is a 
deterministic complex sequence satisfying 
$$
h_{\max}\stackrel{\triangle}{=}\sum_{(k_1,k_2)\in \mathbb{Z}^2} |h(k_1,k_2)| <\infty
$$
and $(U(j_1,j_2),\ (j_1,j_2)\in \mathbb{Z}^2)$ is a sequence of independent
random variables with distribution $\mathcal{CN}(0,1)$.
\end{assump}

\begin{rem} Assumption (A-\ref{hypo1}) is a bit more restrictive than the related assumption 
\cite{BKV96}, which only relies on the summability
of the covariance function of the stationary process.
\end{rem}

For every matrix $A$, we denote by $F^{A\,A^*}$, the empirical
distribution function of the eigenvalues of $A\,A^*$.  Since we will
study at the same time the limiting spectrum of the matrices $Z_n
Z_n^*$ (resp. $(Z_n+A_n)(Z_n+A_n)^*$) and $Z_n^* Z_n$ (resp.
$(Z_n+A_n)^*(Z_n+A_n)$), we can assume without loss of generality that
$c\le 1$. We also assume for simplicity that $N\le n$.

When dealing with vectors, the norm $\|\cdot\|$ will denote the
Euclidean norm. In the case of matrices, the norm $\|\cdot\|$ will
refer to the spectral norm. Denote by $\Cplus$ the set
$\Cplus=\{z\in \C,\ \im(z)>0\}$
and by $C({\mathcal X})$ the set of bounded continuous functions over a given
topological space ${\mathcal X}$ endowed with the supremum norm $\|\cdot \|_{\infty}$.

Let $\mu$ be a probability measure over $\mathbb{R}$. Its Stieltjes transform $f$ is defined by:
$$
f(z)=\int_{\mathbb{R}} \frac{\mu(d\lambda)}{\lambda-z},\quad z\in\Cplus.
$$
We list below the main properties of the Stieltjes transforms that will be needed in the sequel.
\begin{prop} The following properties hold true:
\begin{enumerate}
\item Let $f$ be the Stieltjes transform of $\mu$, then
\begin{itemize}
\item[-] the function $f$ is analytic over $\Cplus$,
\item[-] the function $f$ satisfies: $|f(z)| \le \frac{1}{\im(z)}$,
\item[-] if $z\in\Cplus$ then $f(z)\in \Cplus$,
\item[-] if $\mu(-\infty,0)=0$ then $z\in\Cplus$ implies $z\,f(z)\in \Cplus$.\\
\end{itemize}

\item Conversely, let $f$ be a function analytic over $\mathbb{C}^{+}$
  such that $f(z) \in \mathbb{C}^{+}$ if $z \in \mathbb{C}^{+}$ and
  $|f(z)| |\mbox{Im}(z)|$ bounded on $\mathbb{C}^{+}$. If 
  $\lim_{y\rightarrow +\infty} -iy\,f(iy)=1$, then $f$ is the
  Stieltjes transform of a probability measure $\mu$
  and the following inversion formula holds:
$$
\mu([a,b])=\lim_{\eta\rightarrow 0^+} \frac{1}{\pi} \int_a^b \im f(\xi +\ii \eta)\,d\xi,
$$ 
where $a$ and $b$ are continuity points of $\mu$. If moreover 
  $z f(z) \in \mathbb{C}^{+}$ if $z \in \mathbb{C}^{+}$ then, $\mu(\mathbb{R}^-)=0$.\\

\item Let $\mathbb{P}_n$ and $\mathbb{P}$ be probability measures over $\mathbb{R}$ and denote by $f_n$ and $f$
their Stieltjes transforms. Then 
$$
\left( \forall z\in\Cplus,\ f_n(z) \xrightarrow[n\rightarrow\infty]{} f(z)\right) \quad \Rightarrow \quad 
\mathbb{P}_n\xrightarrow[n\rightarrow\infty]{\mathcal D} \mathbb{P}.  
$$

\end{enumerate}
\end{prop}

Denote by ${\mathcal M}_{\C}({\mathcal X})$ the set of complex measures over the topological set ${\mathcal X}$. 
In the sequel, we will call Stieltjes kernel every application 
$$
\pi: \Cplus  \rightarrow  {\mathcal M}_{\C}({\mathcal X})
$$
either denoted $\pi(z,dx)$ or $\pi_z(dx)$ and satisfying: 
\begin{enumerate}
\item\label{property1} $\forall z\in \Cplus,\ \forall g\in C({\mathcal X})$, 
$$
\left| \int g\,d\pi_z \right|\le \frac{\| g\|_{\infty}}{\im (z)}
$$
\item\label{property2} $\forall g\in C({\mathcal X}),\ \int g\,d\pi_z$ is analytic over $\Cplus$,
\item\label{property3} $\forall z\in \Cplus,\ \forall g\in C({\mathcal X})$ and $g\ge 0$ then 
$\im\left(\int g\,d\pi_z\right)\ge 0$,
\item\label{property4} $\forall z\in \Cplus,\ \forall g\in C({\mathcal X})$ and $g\ge 0$ then 
$\im\left( z\,\int g\,d\pi_z\right) \ge 0$.
\end{enumerate}

\subsection{A quick review of the results for matrices with independent entries}\label{sub-review}
In order to establish the convergence of the empirical distribution of the eigenvalues, we will rely 
on the results based on matrices with independent but not identically distributed entries. Let us recall
here those of interest (the assumptions and the statements are based on \cite{HLN04pre}).

Consider a $N\times n$ random matrix $Y_n$ where the entries are given by 
$$
Y_{j_1 j_2}^{n}=\frac{\Phi(j_1/N,j_2/n)}{\sqrt{n}} X_{j_1 j_2}^{n}
$$
where $X_{j_1 j_2}^n$ and $\Phi$ are defined below.
\begin{assump}\label{hypo2}
The complex random variables $(X_{j_1 j_2}^n\ ;\ 0\le j_1 < N,\,
0\le j_2 < n\,,\, n\ge1)$ are independent and identically distributed (i.i.d.).
 They are centered with
  $\EE | X_{j_1 j_2}^n|^2=1$ and there exists $\epsilon>0$ such that 
  $\EE |X_{j_1 j_2}^n|^{4+\epsilon}<\infty$. 
\end{assump}
\begin{assump}\label{hypo3}
  The function $\Phi: [0,1] \times [0,1] \rightarrow \mathbb{C}$ is such that $|\Phi|^2$ is 
  continuous and therefore there exist a non-negative constant
  $\Phi_{\max}$ such that
\begin{equation}\label{bornitude-var}
\forall (x,y)\in [0,1]^2,\quad 0\le |\Phi(x,y)|^2 \le \Phi^2_{\max}  <\infty.
\end{equation}

\end{assump}
\begin{theo}[independent entries, the centered case \cite{Gir90}]
\label{centered-case-review}
Assume that (A-\ref{hypo2}) and (A-\ref{hypo3}) hold. Then the empirical distribution of the eigenvalues of the matrix 
$Y_n\, Y_n^*$ converges a.s. to a non-random probability measure $\mu$ whose Stieltjes transform $f$ is given by 
$
f(z)=\int_{[0,1]} \pi_z(dx),
$
where $\pi_z$ is the unique Stieljes kernel with support included in $[0,1]$ and satisfying
\begin{equation}\label{stieltjes-centered}
\forall g\in C([0,1]),\quad \int g\,d\pi_z =
\int_0^1 \frac{g(u)}{-z +\int_0^1 \frac{|\Phi|^2(u,t)}{1+c \int_0^1 |\Phi|^2(x,t)\,\pi_z(dx)}dt}du.
\end{equation}
\end{theo}

If one adds a deterministic pseudo-diagonal matrix $\Lambda_n$ to the
matrix $Y_n$, the limiting equation is modified and in fact becomes a
system of equations. 

\begin{assump}\label{hypo4}
  Let $\Lambda_n=(\Lambda_{ij}^n)$ be a complex deterministic $N\times
  n$ matrix whose non-diagonal entries are zero. We assume moreover
  that there exists a probability measure $H(\,du,d\lambda)$ over the
  set $[0,1]\times \mathbb{R}$ with compact support ${\mathcal H}$
  such that
\begin{equation}
\frac1N \sum_{i=1}^N \delta_{\left( \frac iN,\,\left|\Lambda_{ii}^n\right|^2\right)}(du,d\lambda) \xrightarrow[n\rightarrow \infty]{\mathcal D} 
H(\,du,d\lambda).
\end{equation}
\end{assump}
Denote by ${\mathcal H}_c$ the support of the image of probability
measure $H$ under the application $(u,\lambda)\rightarrow
(cu,\lambda)$ and by ${\mathcal R}$ the support of the measure
$\mathbf{1}_{[c,1]}(du) \otimes \delta_0(d\lambda)$ where $\otimes$ denotes the product of measure. The set
$\tilde{\mathcal H}={\mathcal H}_c \cup {\mathcal R}$ will be of
importance in the sequel (see also Remarks 2.4 and 2.5 in
\cite{HLN04pre} for more information).

\begin{theo}[independent entries, the non-centered case \cite{HLN04pre}]
\label{non-centered-case-review}
  Assume that (A-\ref{hypo2}), (A-\ref{hypo3}) and (A-\ref{hypo4})
  hold.  Then the empirical distributions of the eigenvalues of 
  matrices $(Y_n+\Lambda_n)(Y_n+\Lambda_n)^*$ and $(Y_n+\Lambda_n)^*(Y_n+\Lambda_n)$ converge a.s. to 
  non-random probability measures $\mu$ and $\tilde{\mu}$ whose Stieltjes
  transforms $f$ and $\tilde{f}$ are given by
$$
f(z)=\int_{\mathcal H} \pi_z(dx)\qquad \textrm{and}\qquad \tilde{f}(z)=\int_{\mathcal H} \pitilde_z(dx)
$$
where $\pi_z$ and $\pitilde_z$ are the unique Stieljes kernels with supports included in ${\mathcal H}$ and 
$\tilde{\mathcal H}$ and satisfying
\begin{equation}\label{equation1}
\int g\,d\pi_z =\int \frac{g(u,\lambda)}{-z(1+\int |\Phi|^2(u,t)\pitilde(z,dt,d\zeta)) +
\frac {\lambda}{1+c \int |\Phi|^2(t,c u) \pi(z,dt,d\zeta)}} H(du,d\lambda)
\end{equation}
\begin{multline}\label{equation2}
\int g\,d\pitilde_z =c \int \frac{g(c u,\lambda)}{-z(1+c \int |\Phi|^2(t,c u)\pi(z, dt,d\zeta)) +
\frac {\lambda}{1+\int |\Phi|^2(u,t) \pitilde(z, dt,d\zeta)}} H(du,d\lambda)
\\ +(1-c) \int_c^1 \frac{g(u,0)}{-z(1+c \int |\Phi|^2(t,u)\pi(z,dt,d\zeta))} \,du \\
\end{multline}
where (\ref{equation1}) and (\ref{equation2}) hold for every $g\in C({\mathcal H})$
\end{theo}

\section{The limiting distribution in the centered stationary case}
We first introduce the following complex-valued function  
$\Phi:[0,1]\times [0,1] \rightarrow \mathbb{C}$ defined by:
\begin{equation}\label{var-profile}
\Phi(t_1,t_2)=\sum_{(l_1,l_2)\in 
\mathbb{Z}^2} h(l_1,l_2)e^{2\pi\ii (l_1 t_1 - l_2 t_2)} 
\end{equation}
We also introduce the $p\times p$ Fourier matrix 
$F_p=(F_{j_1,j_2}^p)_{0\le j_1,j_2< p}$ defined by: 
\begin{equation}\label{fourier-matrix}
F_{j_1,j_2}^p =
\frac{1}{\sqrt{p}}\exp 2\ii \pi\left( \frac{j_1 j_2}{p}\right).
\end{equation}
Note that matrix $F_p$ is a unitary matrix.
\begin{theo}[stationary entries, the centered case \cite{BKV96,Gir01a}]
  \label{stat-centered} 
  Let $Z_n$ be a $N\times n$ matrix satisfying (A-\ref{hypo1}).  
  Then the empirical
  distribution of the eigenvalues of the matrix $Z_n Z_n^*$ converges
  in probability to the non-random probability measure $\mu$ defined in Theorem \ref{centered-case-review}.
\end{theo}

\subsection{Proof of Theorem \ref{stat-centered}}
Recall that 
$$
Z_{j_1 j_2}^n=
\frac{1}{\sqrt{n}}\sum_{(k_1,k_2)\in \mathbb{Z}^2} 
h(k_1,k_2)U(j_1-k_1,j_2-k_2).
$$
We introduce the $N\times n$ matrix $\tilde{Z}_n$ whose entries are defined by
$$
\tilde{Z}_{j_1 j_2}^n=\frac{1}{\sqrt{n}}\sum_{(k_1,k_2)\in \mathbb{Z}^2} h(k_1,k_2)U(j_1-k_1\ \mathrm{mod}\ N,j_2-k_2 \ \mathrm{mod}\ n).
$$
For simplicity, we shall write $\tilde{U}^n(j_1,j_2)$ instead of 
$U(j_1\ \mathrm{mod}\ N,j_2\ \mathrm{mod}\ n)$. 
Recall that ${\mathcal L}$ stands for the L\'evy distance between
distribution functions.  The main interest in dealing with
matrix $\tilde{Z}_n$ lies in the following two lemmas.

\begin{lemma}\label{lemme:deconv}
Consider the $N\times n$ matrix $Y_n=F_N \tilde{Z}_n F_n^*$. 
Then the entries $Y_{l_1 l_2}^n$ of $Y_n$ can be written 
$$
Y_{l_1 l_2}^n =\frac{1}{\sqrt{n}} 
\Phi\left( \frac{l_1}{N}, \frac{l_2}{n}\right) X^n_{l_1 l_2}
$$
where $\Phi$ is defined in (\ref{var-profile}) and the complex random variables
$\{ X^n_{l_1 l_2}, 0 \leq l_1 < N, 0 \leq l_2 < n \}$ 
are independent with distribution $\mathcal{CN}(0,1)$. 
\end{lemma} 

\begin{proof}[Proof of Lemma \ref{lemme:deconv}]
We first compute the individual entries of matrix 
$Y_n=F_N \tilde{Z}_n F_n^*$:
\begin{eqnarray*}
Y^n_{l_1 l_2}&=&\sum_{\tiny \begin{array}{l}
j_1=0:N-1\\
j_2=0:n-1
\end{array}}
\frac{e^{2\ii \pi\left(\frac{j_1 l_1}{N} - \frac{j_2 l_2}{n}\right)}}
{\sqrt{Nn}} 
\tilde{Z}^n_{j_1 j_2} \\
&=&\frac{1}{\sqrt{n}}\sum_{\tiny \begin{array}{l}
j_1=0:N-1\\
j_2=0:n-1
\end{array}}
\frac{e^{2\ii \pi\left(\frac{j_1 l_1}{N} -
\frac{j_2 l_2}{n}\right)}}{\sqrt{Nn}}
\sum_{(k_1,k_2)\in \mathbb{Z}^2} h(k_1,k_2) \tilde{U}^n(j_1-k_1,j_2-k_2) \\
&=& \frac{1}{\sqrt{n}} \sum_{\tiny \begin{array}{l}
j_1=0:N-1\\
j_2=0:n-1
\end{array}}
\frac{e^{2\ii \pi\left(\frac{j_1 l_1}{N}-\frac{j_2 l_2}{n}\right)}}{\sqrt{Nn}}
\sum_{\tiny \begin{array}{l}
m_1=0:N-1\\
m_2=0:n-1
\end{array}}
U(m_1,m_2) \\
&\phantom{=}& \qquad\times 
\sum_{(k_1,k_2)\in \mathbb{Z}^2} h(j_1-m_1+k_1 N, j_2-m_2 + k_2 n)\\
&=&\frac{1}{\sqrt{n}}\Phi\left(\frac{l_1}{N},\frac{l_2}{n}\right)
\sum_{\tiny \begin{array}{l}
m_1=0:N-1\\
m_2=0:n-1
\end{array}} U(m_1,m_2)  
\frac{e^{2\ii \pi\left(\frac{m_1 l_1}{N}-\frac{m_2 l_2}{n}\right)}}{\sqrt{Nn}}.
\end{eqnarray*}
Let $X_{l_1 l_2}^n$ be the random variable defined as 
$$
X_{l_1 l_2}^n=\sum_{\tiny \begin{array}{l}
m_1=0:N-1\\
m_2=0:n-1
\end{array}} U(m_1,m_2) 
\frac{e^{2\ii \pi\left(\frac{m_1 l_1}{N}-\frac{m_2 l_2}{n}\right)}}
{\sqrt{Nn}} 
$$
for $0\le l_1\le N-1$ and $0\le l_2\le n-1$. Denoting by  
$X_n$ and $U_n$ the $N \times n$ matrices with entries 
$X_{l_1 l_2}^n$ and $U(l_1, l_2)$ respectively, we then have
$X_n = F_N U_n F_n^*$. Define $\mathrm{vec}(A)$ to be the 
vector obtained by stacking the columns of matrix $A$. Then the 
$Nn \times 1$ vectors $\mathbf{X} = \mathrm{vec}(X_n)$ and 
$\mathbf{U} = \mathrm{vec}(U_n)$ are related by the equation
$\mathbf{X} = \left( F_n^* \otimes F_N \right) \mathbf{U}$ (Lemma 4.3.1 in \cite{HorJoh94}), where
$\otimes$ denotes the Kronecker product of matrices. 
The vector $\mathbf{X}$ is a complex Gaussian random vector that
satisfies $\EE \mathbf{X} = 
\left( F_n^* \otimes F_N \right) \ \EE \mathbf{U} = 0$ and 
$\EE \mathbf{X} \mathbf{X}^T = 
\left( F_n^* \otimes F_N \right) \ \EE \mathbf{U} \mathbf{U}^T 
\left( F_n^* \otimes F_N \right) = 0$. After noticing that the matrix
$\left( F_n^* \otimes F_N \right)$ is unitary, we furthermore have
$\EE \mathbf{X} \mathbf{X}^* = 
\left( F_n^* \otimes F_N \right) \ \EE \mathbf{U} \mathbf{U}^*  
\left( F_n^* \otimes F_N \right)^* = I_{nN}$ where $I_{p}$ is the
$p \times p$ identity matrix. In short, the entries of $X_n$ are independent
and have the distribution $\mathcal{CN}(0,1)$. 
Lemma \ref{lemme:deconv} is proved. 
\end{proof}

\begin{lemma}
\label{lemme:approx} 
Let $B_n$ be a $N\times n$ deterministic matrix such that the sequence
$\frac 1n \mathrm{Tr}B_n B_n^*$ is bounded. Then 
  $$
  {\mathcal L}\left( F^{(Z_n +B_n)(Z_n +B_n)^*},F^{(\tilde{Z}_n
      +B_n)(\tilde{Z}_n +B_n)^*}\right) \xrightarrow[n\rightarrow
  \infty]{P} 0,
$$
where $\xrightarrow[]{P}$ denotes convergence in probability.
\end{lemma}

\begin{proof}[Proof of Lemma \ref{lemme:approx}]
Bai's inequality yields:
\begin{multline}\label{ineq:bai}
{\mathcal L}^4(F^{(Z_n+B_n)(Z_n+B_n)^*},F^{(\tilde{Z}_n+B_n)(\tilde{Z}_n+B_n)^*})
\le \frac{2}{n^2} \mathrm{Tr}(Z_n-\tilde{Z_n})(Z_n-\tilde{Z_n})^*\\
\times \mathrm{Tr} 
\left( (Z_n+B_n)(Z_n+B_n)^*+(\tilde{Z}_n+B_n)(\tilde{Z}_n+B_n)^*\right)  
\end{multline}
We introduce the following notations:
\begin{eqnarray*}
\alpha_n &=& \frac 1n \mathrm{Tr}(Z_n -\tilde{Z}_n)(Z_n - \tilde{Z}_n)^*,\\  
\beta_n &=& \frac 1n \mathrm{Tr}(Z_n +B_n)(Z_n +B_n)^*,
\quad \tilde{\beta}_n = \frac 1n \mathrm{Tr}(\tilde{Z}_n +B_n)(\tilde{Z}_n +B_n)^*. 
\end{eqnarray*}
With these notations, Inequality (\ref{ineq:bai}) becomes:
$$
{\mathcal L}^4\left(F^{(Z_n+B_n)(Z_n+B_n)^*},F^{(\tilde{Z}_n+B_n)(\tilde{Z}_n+B_n)^*}\right) \le 2 \alpha_n (\beta_n +\tilde{\beta}_n).
$$
In order to prove that ${\mathcal L}(F^{(Z_n+B_n)(Z_n+B_n)^*}, F^{(\tilde{Z}_n+B_n)(\tilde{Z}_n+B_n)^*})\xrightarrow[]{P} 0$, it is sufficient to prove that 
$\alpha_n (\beta_n +\tilde{\beta}_n)\xrightarrow[]{P}0$, which follows from $\alpha_n\xrightarrow[]{P}0$ and $\beta_n$ 
and $\tilde{\beta}_n$ being tight. Indeed, 
\begin{eqnarray*}
\lefteqn{\mathbb{P}\{ \alpha_n(\beta_n +\tilde{\beta}_n) \ge \epsilon \} 
\le \mathbb{P}\{ \alpha_n \beta_n \ge \epsilon/2\} + 
\mathbb{P}\{ \alpha_n \tilde{\beta}_n \ge \epsilon/2\}} \\
&\le& 
\mathbb{P}\left\{\alpha_n \ge \frac{\epsilon}{2K}\right\}  
+ \mathbb{P}\{\beta_n \ge 2 K\}
+ \mathbb{P}\left\{ \alpha_n \ge \frac{\epsilon}{2\tilde{K}}\right\}
+ \mathbb{P}\{ \tilde{\beta}_n \ge 2 \tilde{K}\}.
\end{eqnarray*}  
Let us first prove that 
\begin{equation}\label{alpha}
\alpha_n\xrightarrow[]{P}0.
\end{equation}
Since $\alpha_n$ is non-negative, it is sufficient by Markov's inequality to prove that $\mathbb{E} \alpha_n \rightarrow 0$.
\begin{eqnarray*}
\alpha_n &=& \frac 1n \mathrm{Tr} (Z_n -\tilde{Z}_n)(Z_n - \tilde{Z}_n)^*\\
&=& \frac{1}{n^2}  \sum_\indicej \left| Z^n_{j_1,j_2} - \tilde{Z}^n_{j_1,j_2}\right|^2\\
&=& \frac{1}{n^2}  \sum_\indicej \Big| \sum_{(k_1,k_2)\in \mathbb{Z}^2} h(k_1,k_2) V(j_1-k_1,j_2-k_2)\Big|^2,
\end{eqnarray*}
where $V(j_1,j_2)$ stands for $U(j_1,j_2)-\tilde{U}^n(j_1,j_2)$. Thus 
$$
\mathbb{E}\,\alpha_n  =\frac{1}{n^2}  \sum_\indicej \sum_{\tiny \begin{array}{l}
(k_1,k_2)\in \mathbb{Z}^2\\
(k'_1,k'_2)\in \mathbb{Z}^2
\end{array}} h(k_1,k_2)h^*(k'_1,k'_2) \mathbb{E} V(j_1-k_1,j_2-k_2)V^*(j_1-k'_1,j_2-k'_2)
$$
Introduce the set ${\mathcal J}=\{0,\cdots,N-1\}\times\{0,\cdots, n-1\}$. Then
\begin{multline*}
\mathbb{E} V(l_1,l_2)V^*(l'_1,l'_2) = \ \mathbf{1}_{\mathbb{Z}^2 -{\mathcal J}}(l_1,l_2) \ \mathbf{1}_{\mathbb{Z}^2 -{\mathcal J}}(l'_1,l'_2)\\
\times \left( \ \mathbf{1}_{(l_1,l_2)}(l'_1,l'_2) + \sum_{(m_1,m_2)\in \mathbb{Z}^2} \ \mathbf{1}_{(l_1,l_2)}(l'_1+m_1 N,l'_2 + m_2 n) \right)
\end{multline*}  
and $\mathbb{E}\, \alpha_n$ becomes $\mathbb{E}\, \alpha_n=\mathbb{E}\, \alpha_{n,1}+\mathbb{E}\, \alpha_{n,2}$ where
\begin{eqnarray*}
\mathbb{E}\, \alpha_{n,1} &=& \frac{1}{n^2} \sum_{\indicej} \sum_{(k_1,k_2)\in \mathbb{Z}^2}  |h(k_1,k_2)|^2 
\ \mathbf{1}_{\mathbb{Z}^2 -{\mathcal J}}(j_1-k_1,j_2-k_2),\\
\mathbb{E}\, \alpha_{n,2} &=& \frac{1}{n^2} \sum_{\indicej} \sum_{\tiny \begin{array}{l}
(k_1,k_2)\in \mathbb{Z}^2\\
(k'_1,k'_2)\in \mathbb{Z}^2
\end{array}} h(k_1,k_2)h^*(k'_1,k'_2) \ \mathbf{1}_{\mathbb{Z}^2 -{\mathcal J}}(j_1-k_1,j_2-k_2)\\
&\phantom{=}&  \times \ \ \mathbf{1}_{\mathbb{Z}^2 -{\mathcal J}}(j_1-k'_1,j_2-k'_2) \sum_{(m_1,m_2)\in \mathbb{Z}^2} 
\ \mathbf{1}_{(k_1,k_2)}(k'_1+m_1 N,k'_2 + m_2 n) 
\end{eqnarray*}
Let us first deal with $\mathbb{E}\,\alpha_{n,2}$. 
\begin{multline*}
\mathbb{E}\, \alpha_{n,2} \le \frac{1}{n^2} \sum_{\indicej} \sum_{(k_1,k_2)\in \mathbb{Z}^2}  |h(k_1,k_2)| 
\ \mathbf{1}_{\mathbb{Z}^2 -{\mathcal J}}(j_1-k_1,j_2-k_2)\\
\times \sum_{(k'_1,k'_2)\in \mathbb{Z}^2}  |h(k'_1,k'_2)|
\ \mathbf{1}_{\mathbb{Z}^2 -{\mathcal J}}(j_1-k'_1,j_2-k'_2)\\
\times \sum_{(m_1,m_2)\in \mathbb{Z}^2} \ \mathbf{1}_{(k_1,k_2)}(k'_1+m_1 N,k'_2 + m_2 n). 
\end{multline*}
Since $h$ is summable over $\mathbb{Z}^2$ by (A-\ref{hypo1}), 
$$
\sum_{(k'_1,k'_2)\in \mathbb{Z}^2}  |h(k'_1,k'_2)|
\ \mathbf{1}_{\mathbb{Z}^2 -{\mathcal J}}(j_1-k'_1,j_2-k'_2)
\sum_{(m_1,m_2)\in \mathbb{Z}^2} \ \mathbf{1}_{(k_1,k_2)}(k'_1+m_1 N,k'_2 + m_2 n)
$$
is bounded by $h_{\max}$ and 
\begin{equation}\label{plug-cesaro}
\mathbb{E}\,\alpha_{n,2} \le \frac{h_{\max}}{n^2}\sum_{\indicej} \sum_{(k_1,k_2)\in \mathbb{Z}^2}  |h(k_1,k_2)| 
\ \mathbf{1}_{\mathbb{Z}^2 -{\mathcal J}}(j_1-k_1,j_2-k_2).
\end{equation}
Since 
$$
\mathbf{1}_{\mathbb{Z}^2 -{\mathcal J}}(j_1-k_1,j_2-k_2)=1 \Leftrightarrow 
\left\{
\begin{array}{ll}
j_1 -k_1 < 0 & \textrm{or}\ j_1 -k_1\ge N,\\
j_2 -k_2 < 0 & \textrm{or}\ j_2 -k_2\ge n
\end{array}\right.
$$
we get:
\begin{eqnarray*}
\lefteqn{\sum_{(k_1,k_2)\in \mathbb{Z}^2}  |h(k_1,k_2)| 
\ \mathbf{1}_{\mathbb{Z}^2 -{\mathcal J}}(j_1-k_1,j_2-k_2)}\\
&=& \sum_{\tiny \begin{array}{l}
k_1=-\infty:j_1-N;\\
k_2=-\infty:j_2-n
\end{array}} |h(k_1,k_2)| 
+\sum_{\tiny \begin{array}{l}
k_1=-\infty:j_1-N;\\
k_2=j_2+1:\infty
\end{array}} |h(k_1,k_2)| \\
&\phantom{=}& +\sum_{\tiny \begin{array}{l}
k_1=j_1+1:\infty;\\
k_2=-\infty:j_2-n
\end{array}} |h(k_1,k_2)| 
+\sum_{\tiny \begin{array}{l}
k_1=j_1+1:\infty;\\
k_2=j_2+1:\infty
\end{array}} |h(k_1,k_2)|. 
\end{eqnarray*}
The changes of variable $
\left\{
\begin{array}{l}
j'_1=N-1-j_1\\
k'_1=-k_1
\end{array}\right.$ and 
$
\left\{
\begin{array}{l}
j'_2=n-1-j_2\\
k'_2=-k_2
\end{array}\right.$ yield
$$
\sum_{\indicej} \sum_{\tiny \begin{array}{l}
k_1=-\infty:j_1-N;\\
k_2=-\infty:j_2-n
\end{array}} |h(k_1,k_2)| 
= \sum_{\tiny \begin{array}{l}
j'_1=0:N-1\\
j'_2=0:n-1
\end{array}} \sum_{\tiny \begin{array}{l}
k'_1=j_1+1:\infty;\\
k'_2=j_2+1:\infty
\end{array}} |h(-k'_1,-k'_2)|.
$$ 
By performing similar changes of variables, one gets:
\begin{eqnarray*}
\lefteqn{\sum_{\indicej} \sum_{(k_1,k_2)\in \mathbb{Z}^2}  |h(k_1,k_2)| 
\ \mathbf{1}_{\mathbb{Z}^2 -{\mathcal J}}(j_1-k_1,j_2-k_2)}\\
&=& \sum_{\indicej} \underbrace{\sum_{\tiny \begin{array}{l}
k_1=j_1+1:\infty;\\
k_2=j_2+1:\infty
\end{array}} |h(-k_1,-k_2)| + |h(-k_1,k_2)| + |h(k_1,-k_2)| + |h(k_1,k_2)|}_{S(j_1,j_2)}.
\end{eqnarray*} 
In order to check that 
\begin{equation}\label{cesaro}
\frac{1}{n^2} \sum_{\indicej} S(j_1,j_2) \xrightarrow[
n\rightarrow\infty\,;\ 
N/n \rightarrow c]{} 0,
\end{equation}
we introduce $T(j)=\sum_{k_1+k_2\ge j+2} |h(-k_1,-k_2)| + |h(-k_1,k_2)| + |h(k_1,-k_2)| +|h( k_1, k_2)| $. Is is straightforward 
to check that $T(j)\xrightarrow[j\rightarrow \infty]{} 0$ and that $S(j_1,j_2)\le T(j_1+j_2)$.
We prove (\ref{cesaro}) by a C\'esaro-like argument: 
Let $n_0\le N$ be such that $T(n_0+1)\le \epsilon$. We have
\begin{equation}\label{cesaro-bis}
\frac{1}{n^2} \sum_{\indicej} S(j_1,j_2) 
=\frac{1}{n^2} \sum_{0\le j_1+j_2\le n_0} S(j_1,j_2) +
\frac{1}{n^2} \sum_{\tiny \begin{array}{l}
n_0+1\le j_1 +j_2;\\
j_1\le N-1,\ j_2\le n-1
\end{array}} S(j_1,j_2).   
\end{equation} 
If $n$ is large enough, then the first part of the right handside of (\ref{cesaro-bis}) is lower than $\epsilon$. Moreover, 
$$
\frac{1}{n^2} \sum_{\tiny \begin{array}{l}
n_0+1\le j_1 +j_2;\\
j_1\le N-1,\ j_2\le n-1
\end{array}} S(j_1,j_2)\le \frac{1}{n^2} \sum_{\tiny \begin{array}{l}
n_0+1\le j_1 +j_2;\\
j_1\le N-1,\ j_2\le n-1
\end{array}} T(n_0+1)\le \epsilon
$$
and (\ref{cesaro}) is proved. By pluging (\ref{cesaro}) into (\ref{plug-cesaro}), we prove that $\mathbb{E}\, \alpha_{n,2} \rightarrow 0$.
Using the same kind of arguments, one proves that $\mathbb{E}\, \alpha_{n,1} \rightarrow 0$. Finally, (\ref{alpha}) is proved:
$\alpha_n \xrightarrow[]{P} 0$.
\\ 
Let us now check that 
\begin{equation}
\label{tightness}
\exists\, K>0,\quad \mathbb{E}\, \beta_n \le K\qquad \textrm{and}\quad \exists\, \tilde{K}>0,\quad \mathbb{E}\, 
\tilde{\beta}_n \le \tilde{K}.
\end{equation} 
This will imply the tightness of $\beta_n$ and $\tilde{\beta}_n$. 
\\ 
Recall that by assumption there exists $B_{\max}$ such that 
$\sup_n \frac 1n \mathrm{Tr}B_n B_n^* \le B_{\max}$. Consider now:
\begin{eqnarray*}
\frac 1n
\mathrm{Tr}(Z_n + B_n)(Z_n +B_n)^* &\le& 
\left( \left(\frac 1n \mathrm{Tr} Z_n\,Z_n^*\right)^{\frac 12}
+\left(\frac 1n \mathrm{Tr}B_n\,B_n^*\right)^{\frac 12} \right)^2 
\\
&\le & 
\left( \left(\frac 1n \mathrm{Tr} Z_n\,Z_n^*\right)^{\frac 12}
  +B_{\max}^{\frac 12} \right)^2
\end{eqnarray*}
In particular,
\begin{eqnarray}
\mathbb{E} \frac{\mathrm{Tr}(Z_n + B_n)(Z_n +B_n)^*}{n} &\le &
\mathbb{E} \frac{\mathrm{Tr} Z_n Z_n^*}{n} 
+2B_{\max}^{\frac 12}\mathbb{E}\left( \frac{\mathrm{Tr} Z_n Z_n^*}{n}\right)^{\frac 12} + B_{\max}\nonumber \\
&\stackrel{(a)}{\le}& 
\mathbb{E} \frac{\mathrm{Tr} Z_n Z_n^*}{n}
+2B_{\max}^{\frac 12} \left(\mathbb{E}\left( \frac{\mathrm{Tr} Z_n Z_n^*}{n}\right)\right)^{\frac 12} + B_{\max}\label{controle-trace}
\end{eqnarray}
where $(a)$ follows from Jensen's inequality. Notice that (\ref{controle-trace}) still holds if one replaces $Z_n$ by $\tilde{Z}_n$.
Therefore in order to prove (\ref{tightness}), it is sufficient to prove that:
$$
\exists \ K'>0,\quad \mathbb{E}\left( \frac{\mathrm{Tr} Z_n Z_n^*}{n}\right)\le K'\quad \textrm{and}\quad 
\exists \ \tilde{K}'>0,\quad \mathbb{E}\left( \frac{\mathrm{Tr} \tilde{Z}_n \tilde{Z}_n^*}{n}\right)\le \tilde{K}'.
$$
Consider 
$$
\mathbb{E}\left( \frac{\mathrm{Tr} Z_n Z_n^*}{n}\right)= \frac 1n \sum_{\tiny \begin{array}{l}
j_1=1:N\\
j_2=1:n
\end{array}} \mathbb{E}|Z^n_{j_1 j_2}|^2 =N \mathbb{E}|Z^n_{1 1}|^2=\frac Nn C(0,0),
$$
where $C$ is defined by (\ref{variance}). 
This quantity is asymptotically bounded. From lemma \ref{lemme:deconv}, we 
have
$$
\mathbb{E}\left( \frac{\mathrm{Tr} \tilde{Z}_n \tilde{Z}_n^*}{n}\right)= 
\mathbb{E}\left( \frac{\mathrm{Tr} Y_n Y_n^*}{n}\right)=
\frac{1}{n^2} \sum_{\tiny \begin{array}{l}
j_1=1:N\\
j_2=1:n
\end{array}} \left|\Phi\left(\frac{j_1}N,\frac{j_2}n\right)\right|^2 \mathbb{E}|X^n_{j_1 j_2}|^2 \le \frac Nn \Phi_{\max}^2,
$$
which is also asymptotically bounded. Eq. (\ref{tightness}) is proved and so is Lemma \ref{lemme:approx}.
\end{proof}

\begin{proof}[Proof of Theorem \ref{stat-centered}]
Lemma \ref{lemme:approx} implies that 
\begin{equation}\label{proba1-0}
\mathbb{P}\left\{ {\mathcal L}\left(F^{Z_n Z_n^*},F^{\tilde{Z}_n \tilde{Z}_n^*}\right) \ge \epsilon\right\}
\xrightarrow[n\rightarrow\infty]{}0\quad \textrm{for\ every}\ \epsilon>0. 
\end{equation}
By lemma \ref{lemme:deconv},
$
F_N \tilde{Z}_n \tilde{Z}_n^* F_N^* 
= Y_n Y_n^* .$
Since $F_N$ is unitary, $\tilde{Z}_n \tilde{Z}_n^*$ and
$Y_n Y_n^*$ have the same eigenvalues. Moreover, matrix $Y_n$ 
fulfills (A-\ref{hypo2}) and the variance profile $\Phi$ defined in (\ref{var-profile}) satisfies (A-\ref{hypo3}) 
since $(h(k_1,k_2)\in (k_1,k_2)\mathbb{Z}^2)$ is summable; therefore one can apply Theorem 
\ref{centered-case-review}. In particular,
\begin{equation}\label{proba2-0}
F^{\tilde{Z}_n \tilde{Z}_n^*} \xrightarrow[n\rightarrow \infty]{} \mu\quad a.s. \qquad \Longrightarrow \qquad 
\forall \epsilon>0,\quad  \mathbb{P}\left\{ {\mathcal L}\left(F^{\tilde{Z}_n \tilde{Z}_n^*},\mu\right) 
\ge \epsilon\right\}
\xrightarrow[n\rightarrow\infty]{}0
\end{equation}
where $\mu$ is the probability distribution defined in Theorem \ref{centered-case-review}.
Eq. (\ref{proba1-0}) together with (\ref{proba2-0}) imply that $F^{Z_n Z_n^*} \xrightarrow[]{P} \mu$ and
Theorem \ref{stat-centered} is proved.
\end{proof}

\section{The limiting distribution in the non-centered stationary case}
Recall the definitions of function $\Phi$  and matrix $F_p$ (respectively 
defined in (\ref{var-profile}) and (\ref{fourier-matrix})).
\begin{theo}[stationary entries, the non-centered case]
\label{stat-noncentered}
  Let $Z_n$ be a $N\times n$ matrix satisfying (A-\ref{hypo1}); let
  $A_n$ be a $N \times n$ matrix such that $\Lambda_n  = F_N A_n F_n^*$ is 
  $N\times n$ pseudo-diagonal and satisfies (A-\ref{hypo4}). 
  Then the empirical distributions of the eigenvalues
  of matrices $(Z_n+A_n)(Z_n+A_n)^*$ and $(Z_n+A_n)^*(Z_n+A_n)$ converge in probability to the
  non-random probability measures $\mu$ and $\tilde{\mu}$ defined in 
  Theorem \ref{non-centered-case-review}.
\end{theo}

\begin{proof}[Proof of Theorem \ref{stat-noncentered}]
We denote by $F^n=F^{(Z_n+A_n)(Z_n+A_n)^*}$ and $\tilde{F}^n=F^{(\tilde{Z}_n+A_n)(\tilde{Z}_n+A_n)^*}$. 
Since $\Lambda_n$ satisfies (A-\ref{hypo4}), $\frac 1n \mathrm{Tr}A_n A_n^*=\frac 1n \mathrm{Tr}\Lambda_n \Lambda_n^*$
is bounded and Lemma \ref{lemme:approx} implies that 
\begin{equation}\label{proba1}
\mathbb{P}\left\{ |{\mathcal L}(F^n,\tilde{F}^n)|\ge \epsilon\right\}
\xrightarrow[n\rightarrow\infty]{}0\quad \textrm{for\ every}\ \epsilon>0. 
\end{equation}
By lemma \ref{lemme:deconv} and the assumption over $A_n$, 
$$
(\tilde{Z}_n+A_n)(\tilde{Z}_n+A_n)^*=F_N (Y_n +\Lambda_n)(Y_n +\Lambda_n)^* F_N^*.
$$
Since the Fourier matrix $F_N$ is unitary, $(\tilde{Z}_n+A_n)(\tilde{Z}_n+A_n)^*$ and
$(Y_n +\Lambda_n)(Y_n +\Lambda_n)^*$ have the same eigenvalues. 
Since $\Phi$ defined in (\ref{var-profile}) satisfies (A-\ref{hypo3}), 
the matrices $Y_n$ and $\Lambda_n$ 
fulfill assumptions (A-\ref{hypo2}), (A-\ref{hypo3}) and (A-\ref{hypo4}) therefore one can apply Theorem 
\ref{non-centered-case-review}. In particular,
\begin{equation}\label{proba2}
\tilde{F}^n \xrightarrow[n\rightarrow \infty]{} \mu\quad a.s. \qquad \Longrightarrow \qquad 
\forall \epsilon>0,\quad  \mathbb{P}\left\{ |{\mathcal L}(\tilde{F}^n,\mu)|\ge \epsilon\right\}
\xrightarrow[n\rightarrow\infty]{}0
\end{equation}
where $\mu$ is the probability distribution defined in Theorem \ref{non-centered-case-review}.
Eq. (\ref{proba1}) together with (\ref{proba2}) imply that $F^n \xrightarrow[]{\mu} \mathbb{P}$ and
Theorem \ref{stat-noncentered} is proved.
\end{proof}

In the square case $n\times n$, we can deal with slightly more general matrices $A_n$.
\begin{assump}
\label{hypo5}
The $n\times n$ matrix $A_n$ is a Toeplitz matrix
defined as $A_n=(a(j_1-j_2))_{0\le j_1,j_2 < n}$ where 
$(a(j))_{j\in \mathbb{Z}}$ is a deterministic sequence of complex numbers
satisfying:
$$
\sum_{j\in \mathbb{Z}} |a(j)| <\infty.
$$ 
\end{assump}
Let $\psi:[0,1]\mapsto \mathbb{C}$ be the so called symbol of $A_n$ 
defined as 
\begin{equation}
\label{eq:def-psi} 
\psi(t)=\sum_{j \in \mathbb{Z}} a(j)e^{2 \ii \pi j t}.
\end{equation} 
Due to (A-\ref{hypo5}), $\psi$ is bounded and continuous.
\begin{theo}[stationary entries, the non-centered square case]
\label{square}
  Let $Z_n$ be a $n\times n$ matrix satisfying (A-\ref{hypo1}); let
  $A_n$ be a $n\times n$ matrix satisfying (A-\ref{hypo5}). 
  Then the empirical distributions of the eigenvalues
  of matrices $(Z_n+A_n)(Z_n+A_n)^*$ and $(Z_n+A_n)^*(Z_n+A_n)$ converge in probability to 
  non-random probability measures $\mu$ and $\tilde{\mu}$ whose Stieltjes
  transforms $f$ and $\tilde{f}$ are  given by 
$$
f(z)=\int_{[0,1]} \pi_z(dx)\qquad \textrm{and}\qquad \tilde{f}(z)=\int_{[0,1]} \pitilde_z(dx)
$$
  where $\pi_z$ and $\pitilde_z$ are the unique Stieltjes kernels with 
  supports included in $[0,1]$ and satisfying
  the system of equations:
\begin{eqnarray}
\int g\,d\pi_z &=&
\int_0^1 \frac{g(u)}{-z(1+\int |\Phi(u,\cdot)|^2 d\pitilde_z) +
\frac {|\psi(u)|^2}{1+ \int |\Phi(\cdot, u)|^2 d\pi_z}}du
\label{ker1} \\
\int g\,d\pitilde_z &=& 
\int_0^1 \frac{g(u)}{-z(1+ \int |\Phi(\cdot , u)|^2d\pi_z) +
\frac {|\psi(u)|^2}{1+\int |\Phi(u,\cdot)|^2 d\pitilde_z}} du
\label{ker2}
\end{eqnarray}
for every function $g\in C([0,1])$.
\end{theo}

\begin{proof}
The proof is based on the fact that a Toeplitz matrix $A_n$ is very close 
to a Toeplitz circulant matrix $\tilde{A}_n$ defined in such a way that the 
diagonal matrix $\Lambda_n = F_n \tilde{A}_n F_n^*$ 
satisfies assumption (A-4). Denoting by $\psi_n$ the truncated function 
$\psi_n(t)=\sum_{j=-n}^n a(j) e^{2 \ii \pi j t }$, we choose 
$\tilde{A}_n$ to be the matrix whose entries are defined by 
$$
\tilde{a}^n_{j_1 j_2} 
=\frac 1n 
\sum_{k=0}^{n-1} \psi_n\left(\frac{k}{n} \right) 
\exp\left(\frac{- 2\pi \ii k(j_1 -j_2)}{n}\right)
.
$$ 
Notice that in this case, $\Lambda_n = F_n \tilde{A}_n F_n^*$ is given
by $\Lambda_n = \mathrm{diag}\left( \left[ \psi_n(0), \psi_n(\frac 1n), \ldots,
\psi_n(\frac{n-1}{n}) \right] \right)$ where $\mathrm{diag}(v)$ is the diagonal
matrix bearing the entries of the vector $v$ on its diagonal. \\
One can also prove that the complex number
$\tilde{a}^n(j_1 - j_2) = \tilde{a}^n_{j_1 j_2}$ satisfies 
$\tilde{a}^n(0)=a(0)+a(n)+a(-n)$ and
$$
\tilde{a}^n(j)
=\left\{
\begin{array}{ll}
a(j)+a(j-n) &\textrm{if}\ n-1\ge j>0,\\
a(j)+a(j+n) &\textrm{if}\ -n+1\le j<0.
\end{array}\right.
$$
We denote by $F^n$ and $\breve{F}^n$ the distribution
functions $F^n=F^{(Z_n+A_n)(Z_n+A_n)^*}$ and 
$\breve{F}^n=F^{(Z_n+\tilde{A}_n)(Z_n+\tilde{A}_n)^*}$. 
We shall prove that ${\mathcal L}(F^n,\breve{F}^n) \to 0$ as
$n \to \infty$. \\
Bai's inequality yields: 
\begin{equation}
\label{controle0}
{\mathcal L}^4(F^n,\breve{F}^n)
\le 
\frac{2}{n^2}\mathrm{Tr}(A_n -\tilde{A}_n)(A_n -\tilde{A}_n)^* 
\mathrm{Tr}(A_n A_n^* + \tilde{A}_n \tilde{A}_n^*). 
\end{equation}
We first prove that $n^{-1}\mathrm{Tr}(A_n A_n^*)$ and 
$n^{-1}\mathrm{Tr}(\tilde{A}_n \tilde{A}_n^*)$ are bounded:
\begin{equation}\label{controle1}
\frac 1n \mathrm{Tr}A_n A_n^* = \frac 1n \sum_{j_1,j_2=0}^{n-1} |a(j_1-j_2)|^2 
=\sum_{j=-n+1}^{n-1} |a(j)|^2\left(1 -\frac{|j|}n\right)\le 
\Big(\sum_{j\in \mathbb{Z}} |a(j)|\Big)^2.
\end{equation}
Moreover, 
\begin{equation}
\label{controle2}
\frac 1n \mathrm{Tr}\tilde{A}_n \tilde{A}_n^* =
\frac 1n \mathrm{Tr}\Lambda_n \Lambda_n^* =
\frac 1n 
\sum_{j=0}^{n-1} \left|\psi_n\left(\frac jn\right)\right|^2
\le \Big(\sum_{j\in \mathbb{Z}} |a(j)|\Big)^2.
\end{equation}
We now prove that 
\begin{equation}\label{controle3}
\frac 1n \mathrm{Tr} (A_n -\tilde{A}_n)(A_n -\tilde{A}_n)^* 
\xrightarrow[n\rightarrow\infty]{} 0.
\end{equation}
Indeed, 
\begin{eqnarray*}
\lefteqn{\frac 1n \mathrm{Tr} (A_n -\tilde{A}_n)(A_n -\tilde{A}_n)^*
= \frac 1n \sum_{j_1,j_2=0}^{n-1} |a(j_1 -j_2)-\tilde{a}^n(j_1-j_2)|^2}\\
&=& \sum_{j=-(n-1)}^{n-1} |a(j)-\tilde{a}^n(j)|^2\left(1 -\frac{|j|}n\right)\\
&=& |a(-n)+a(n)|^2 + \sum_{j=1}^{n-1} \Big(|a(j-n)|^2+|a(n-j)|^2\Big)\left(1 -\frac{j}n\right) \\
&=& |a(-n)+a(n)|^2 + \sum_{j=1}^{n-1}\frac{j}n \Big(|a(j)|^2+|a(-j)|^2\Big) \\
&\le & |a(-n)+a(n)|^2 +\frac 1n \sum_{j=1}^J j \Big(|a(j)|^2+|a(-j)|^2\Big)
+  \sum_{j=J+1}^\infty \Big(|a(j)|^2+|a(-j)|^2\Big)
\end{eqnarray*}
By first taking $J$ large enough then $n$ large enough, the claim is proved 
by a $2\epsilon$-argument. \\ 
Eq. (\ref{controle0}) together with the arguments provided by
(\ref{controle1}), (\ref{controle2}) and (\ref{controle3}) imply that
$${\mathcal L}(F^n,\breve{F}^n) \xrightarrow[n \to \infty]{ } 0.$$
It remains to prove that $\breve{F}^n$ converges towards the non random 
probability distribution characterized by equations 
(\ref{ker1}) and (\ref{ker2}).  As previously, the variance profile $\Phi$  
defined in (\ref{var-profile}) satisfies (A-\ref{hypo3}). Moreover, we have 
$$
\frac 1n \sum_{i=1}^n 
\delta_{\left(\frac in, 
\left|\psi_n\left(\frac{i-1}n\right)\right|^2\right)}
\xrightarrow[n\rightarrow\infty]{ } H(du,d\lambda) 
$$
where $H(du,d\lambda)$ is the image of the Lebesgue measure over $[0,1]$
under $u\mapsto (u, |\psi(u)|^2)$. 
Therefore $\Lambda_n$ satisfies (A-4) and 
Theorem \ref{stat-noncentered} can be applied. This completes the proof of
Theorem \ref{square}. 
\end{proof}

\section{Remarks on the real case}\label{section:real}
In the case where the entries of matrix $Z_n$ are given by 
$$
Z^n_{j_1 j_2}=\frac{1}{\sqrt{n}}\sum_{(k_1,k_2)\in \mathbb{Z}^2} h(k_1,k_2)U(j_1-k_1,j_2-k_2),
$$ 
where $(h(k_1,k_2),\ (k_1,k_2)\in \ZZ^2)$ is a deterministic real and summable 
sequence and where $U(j_1,j_2)$ are real standard independent gaussian r.v.'s, the conclusion of Lemma \ref{lemme:deconv}
is no longer valid. In fact the entries of $Y_n=F_N \tilde{Z}_n F_n^*$ are far from being independent since
straightforward computation yields:
$$
Y^n_{l_1, l_2} = Y^{n^*}_{N-l_1, n-l_2}\quad \textrm{for}\quad  0< l_1 < N\ \textrm{and}\ 0 < l_2 < n.
$$
We introduce the $p\times p$ orthogonal matrix $Q_p=(Q^P_{j_1 j_2})_{0\le j_1,j_2<p}$ defined by:
$$
Q^p_{0,j_2}=\frac{1}{\sqrt{p}},\quad 0\le j_2<p.
$$
In the case where $p$ is even, the entries $Q^p(j_1,j_2)$ ($j_1\ge 1$) are defined by  
$$
\left\{
\begin{array}{ll}
Q^p_{2 j_1 -1,j_2} =\sqrt{\frac 2p} \cos\left(\frac{2\pi j_1 j_2}{p}\right) & \textrm{if}\quad 1\le j_1\le \frac p2 -1, 0\le j_2<p;\\
Q^p_{2 j_1,j_2} =\sqrt{\frac 2p} \sin\left(\frac{2\pi j_1 j_2}{p}\right) & \textrm{if}\quad 1\le j_1\le \frac p2 -1, 0\le j_2<p;\\
Q^p_{p-1,j_2}= \frac{(-1)^{j_2}}{\sqrt{p}} & \textrm{if}\quad 0\le j_2<p.
\end{array}\right.
$$
In the case where $p$ is odd, they are defined by 
$$
\left\{
\begin{array}{ll}
Q^p_{2 j_1 -1,j_2} =\sqrt{\frac 2p} \cos\left(\frac{2\pi j_1 j_2}{p}\right) & \textrm{if}\quad 1\le j_1\le \frac{p-1}2, 0\le j_2<p;\\
Q^p_{2 j_1,j_2} =\sqrt{\frac 2p} \sin\left(\frac{2\pi j_1 j_2}{p}\right) & \textrm{if}\quad 1\le j_1\le \frac{p-1}2, 0\le j_2<p.
\end{array}\right.
$$
In the sequel, $\lfloor x\rfloor$ stands for the integer part of
$x$. The following result is the counterpart of Lemma
\ref{lemme:deconv} in the real case.
\begin{lemma}
\label{lemme:deconv-reel}
Consider the $N\times n$ matrix $W_n= Q_N \tilde{Z}_n Q_n^\T$ where
$A^\T$ is the transpose of matrix $A$. 
Then the entries $W_{l_1 l_2}^n$ of $W_n$ can be written as 
$$
W_{l_1 l_2}^n =\frac{1}{\sqrt{n}} 
\left| \Phi\left( \frac{1}{N} \left\lfloor \frac{l_1+1}{2} \right\rfloor, 
\frac{1}{n} \left\lfloor \frac{l_2+1}{2} \right\rfloor \right) \right| 
X^n_{l_1 l_2}
$$
where $\Phi$ is defined in (\ref{var-profile}) and the real random variables
$\{ X^n_{l_1 l_2}, 0 \leq l_1 < N, 0 \leq l_2 < n \}$ 
are independent standard gaussian r.v.'s.
\end{lemma} 
The proof is computationally more involved but similar in spirit to that of Lemma \ref{lemme:deconv}. It is thus ommited.

As a consequence of this lemma, Theorems \ref{stat-centered} and 
\ref{stat-noncentered} remain true with the following minor modification:  
In Eq. (\ref{stieltjes-centered}), (\ref{equation1}) and (\ref{equation2}), 
the quantity $| \Phi |^2$ must be replaced by $\Phi_{\mathrm{R}}^2$ where
$$
\Phi_{\mathrm{R}}(u, v) = \left| \Phi(u/2, v/2) \right|.
$$ 
Similarly, in the case where the Toeplitz matrix $A_n$ introduced in (A-\ref{hypo5}) is real, 
Theorem \ref{square} remains true if one replaces in (\ref{ker1}) and (\ref{ker2}) 
the quantities $|\Phi|^2$ and $|\psi|^2$ by $\Phi_{\mathrm{R}}^2$ and $\psi^2_{\mathrm{R}}$ where
$$
\psi_{\mathrm{R}}(u) = | \psi(u/2) |.
$$
The proof of Theorem \ref{square} can be modified by replacing
the Fourier matrices $F_p$ by $Q_p$ (see also \cite{BroDav91}, chap. 4 for elements about 
the pseudo-diagonalization of a real Toeplitz matrix via real orthogonal matrices $Q_p$).

\nocite{Bai99}

\bibliography{math}

\noindent {\sc Walid Hachem},\\
Sup\'elec (Ecole Sup\'erieure d'Electricit\'e)\\
Plateau de Moulon, 3 rue Joliot-Curie\\
91192  Gif Sur Yvette Cedex, France.\\
e-mail: walid.hachem@supelec.fr\\
\\
\noindent {\sc Philippe Loubaton},\\
IGM LabInfo, UMR 8049, Institut Gaspard Monge,\\
Universit\'e de Marne La Vall\'ee, France.\\
5, Bd Descartes, Champs sur Marne, \\
77454 Marne La Vall\'ee Cedex 2, France.\\
e-mail: loubaton@univ-mlv.fr\\
\\
\noindent {\sc Jamal Najim},\\ 
CNRS, T\'el\'ecom Paris\\ 
46, rue Barrault, 75013 Paris, France.\\
e-mail: najim@tsi.enst.fr

\end{document}